\documentclass[a4paper,11pt]{amsart}

\usepackage{graphicx}
\usepackage{mathptmx}
\usepackage{amsmath}
\usepackage{amssymb}
\usepackage{enumitem}
\usepackage{xcolor}

\newmuskip\pFqmuskip

\newcommand*\pFq[6][8]{%
  \begingroup 
  \pFqmuskip=#1mu\relax
  \mathcode`=\string"8000
  \begingroup\lccode`\~=`\,
  \lowercase{\endgroup\let~}\pFqcomma
  F^{#2}_{#3}{\left(\genfrac..{0pt}{}{#4}{#5}\bigg|#6\right)}%
  \endgroup
}
\newcommand{\pFqcomma}{\mskip\pFqmuskip}

\newtheorem{theorem}{Theorem}[section]

\begin{document}

\title[Probabilistic degenerate Bernstein polynomials]{Probabilistic degenerate Bernstein polynomials}

\author{Jinyu Wang}
\address{School of Science, Xi'an Technological University, Xi'an, 710021, Shaanxi, China}
\email{2356246083@qq.com}

\author{Yuankui Ma*}
\address{School of Science, Xi'an Technological University, Xi'an, 710021, Shaanxi, China}
\email{mayuankui@xatu.edu.cn}

\author{Taekyun  Kim*}
\address{School of Science, Xi'an Technological University, Xi'an, 710021, Shaanxi, China; Department of Mathematics, Kwangwoon University, Seoul 139-701, Republic of Korea}
\email{tkkim@kw.ac.kr}
\author{Dae San  Kim*}
\address{Department of Mathematics, Sogang University, Seoul 121-742, Republic of Korea}
\email{dskim@sogang.ac.kr}

\thanks{ * corresponding authors}

\begin{abstract}
In recent years, both degenerate versions and probabilistic extensions of many special numbers and polynomials have been explored. For instance, degenerate Bernstein polynomials and probabilistic Bernstein polynomials were investigated earlier. Assume that $Y$ is a random variable whose moment generating function exists in a neighborhood of the origin. The aim of this paper is to study probabilistic degenerate Bernstein polynomials associated with $Y$ which are both probabilistic extension of the degenerate Bernstein polynomials and degenerate version of the probabilistic Bernstein polynomials associated with $Y$. We derive several explicit expressions and certain related identities for those polynomials. In addition, we treat the special cases of the Poisson random variable, the Bernoulli random variable and of the binomial random variable.
\end{abstract}

\maketitle

\section{Introduction}
In recent years, various degenerate versions of many special numbers and polynomials have been studied and a lot of interesting results have been discovered (see [6,8,10,15,16,19-21,24] and the references therein), since Carlitz initiated a study of degenerate versions of Bernoulli and Euler polynomials, namely degenerate Bernoulli and degenerate Euler numbers (see [3]). These degenerate versions have been explored by means of generating functions, combinatorial methods, umbral calculus, $p$-adic analysis, probability theory, special functions and differential equations. Furthermore, probabilistic extensions of many special numbers and polynomials have been investigated as well (see [4,7,13,14,17,23,26,27] and the references therein). \par
Bernstein polynomials were used by Bernstein to give a constructive proof for the
Stone–Weierstrass approximation theorem in 1912. About half a century later, they were used to design automobile bodies at Renault by Pierre B\'{e}zier. Also, the Bernstein polynomials are the mathematical basis for B\'{e}zier curves, which are frequently used in computer graphics and related fields such as animation, modeling, CAD, and CAGD. \par
Recently, Karagenc et al studied the probabilistic Bernstein polynomials (see [7]) which are probabilistic extensions of the Bernstein polynomials. The aim of this paper is to study the  probabilistic degenerate Bernstein polynomials (see \eqref{18}). These are both degenerate version of the probabilistic Bernstein polynomials (see [7]) and probabilistic extension of the degenerate Bernstein polynomials (see [15,16,20]). We derive several explicit expressions and certain  related identities for $B_{k,n}^{Y}(x|\lambda)$. We also treat the special cases of the Poisson random variable with parameter $\alpha (>0)$, the Bernoulli random variable with probability of success $p$, and of the binomial random variable with parameters $m, p(>0)$. \par
The outline of this paper is as follows. In Section 1, we recall the Bernstein polynomials which are the probabilistic mass function of the binomial distribution. We remind the reader of the Stirling numbers of the first kind $S_{1}(n,k)$ and of the second kind $S_{2}(n,k)$. We recall the degenerate exponentials, the degenerate Stirling numbers of the second kind ${ n \brace k}_{\lambda}$, the degenerate Bell polynomials $\phi_{n,\lambda}(x)$ and the degenerate Bernstein polynomials $B_{k,n}(x|\lambda)$. Let $Y$ be a random variable such that the moment generating function of $Y$ exists in a neighborhood of the origin (see \eqref{7-1}). We remind the reader of the probabilistic degenerate Stirling numbers of the second kind ${n \brace k}_{Y,\lambda}$ associated with $Y$, the probabilistic degenerate Bell polynomials $\phi_{n,\lambda}^{Y}(x)$ associated with $Y$,  the probabilistic degenerate Bernoulli polynomials associated $\beta_{n,\lambda}^{Y}(x)$ with $Y$ and the probabilistic degenerate Euler polynomials $\mathcal{E}_{n,\lambda}^{Y}(x)$ associated with $Y$. Finally, we recall the probabilistic Bernstein polynomials introduced by Karagenc et al. Section 2 contains the main results of this paper. We first define the probabilistic degenerate Bernstein polynomials $B_{k,n}^{Y}(x|\lambda)$ associated with $Y$. In Theorem 2.1, we derive an explicit expression for $B_{k,n}^{Y}(x|\lambda)$ as a sum involving ${n-k \brace j}_{Y,\lambda}$. In Theorem 2.2, we get an expression for $B_{k,n}^{Y}(x|\lambda)$ as a sum involving $\beta_{n-m,\lambda}^{(k,Y)}(1-x)$. We obtain an explicit expression for $B_{k,n+k}^{Y}(x|\lambda)$ as a sum involving ${n \brace j}_{Y,\lambda}$ in Theorem 2.3. In Theorem 2.4, the expectation $E\big[(Y)_{n-k,\lambda}\big]$
is expressed as a sum involving $B_{k,n-l}^{Y}(x|\lambda)$. We represent $\mathcal{E}_{n-k,\lambda}^{Y}(1-x)$ as a sum involving $\mathcal{E}_{n-l,\lambda}^{Y}$ and $B_{k,l}^{Y}(x|\lambda)$ in Theorem 2.5. We express $B_{k,n}^{Y}(x|\lambda)$ as a sum involving $\beta_{n-m,\lambda}^{(k,Y)}(kx)$ and ${m \brace k+1}_{Y,\lambda}$ in Theorem 2.6. In Theorem 2.7, $B_{k,n}^{Y}(x|\lambda)$ is represented as a sum involving $S_{1}(j,l)$ and ${n-k \brace j}_{Y,\lambda}$. When $Y$ is the Poisson random variable with parameter $\alpha (>0)$, $B_{k,n}^{Y}(x|\lambda)$ is expressed in terms of $\phi_{n-k,\lambda}\big(\alpha (1-x)\big)$ in Theorem 2.8 and of ${n-k \brace j}_{\lambda}$ in Theorem 2.9. $B_{k,n}^{Y}(x|\lambda)$ are expressed in terms of the degenerate Stirling numbers of the second kind when $Y$ is the Bernoulli random variable with probability of success $p$ (see Theorem 2.10) and when $Y$ is the binomial random variable with parameters $m,\,p$ (see Theorem 2.11). For the rest of this section, we recall the facts that are needed throughout this paper. \par
\vspace{0.1in}
It is well known that the Bernstein polynomials of degree $n$ are defined by
\begin{displaymath}
B_{k,n}(x)=\binom{n}{k}x^{k}(1-x)^{n-k},\quad (n\ge k\ge 0),\quad (\mathrm{see}\ [7,16,20]).
\end{displaymath}
A Bernoulli trial involves performing a random experiment and noting whether a particular event $A$ occurs. The outcome of Bernoulli trial is said to be ``success" if $A$ occurs and a ``failure" otherwise. For $n=1,2,\dots$ and $0 \le p \le 1$, the probability $P_{n}(k)$ of $k$ success in $n$ independent Bernoulli trials is given by the binomial probability law :
\begin{displaymath}
	P_{n}(k)=\binom{n}{k}p^{k}(1-p)^{n-k},\quad\mathrm{for}\quad k=0,1,2,\dots,n,\quad (\mathrm{see}\ [7,16,20]).
\end{displaymath}
We note that the Bernstein polynomials are the probabilistic mass function of the binomial distribution with parameters $(n,\ x=p)$.\par
For any $\lambda\in\mathbb{R}$, the degenerate exponentials are defined by
\begin{equation}
e_{\lambda}^{x}(t)=\sum_{k=0}^{\infty}(x)_{k,\lambda}\frac{t^{n}}{n!},\quad e_{\lambda}(t)=e_{\lambda}^{1}(t),\quad (\mathrm{see}\ [8-20]), \label{1}
\end{equation}
where $(x)_{0,\lambda}=1,\ (x)_{n,\lambda}=x(x-\lambda)(x-2\lambda)\big(x-(k-1)\lambda\big),\ (k\ge 1)$. \par
For $n\ge 0$, the Stirling numbers of the first kind are given by
\begin{equation}
(x)_{n}=\sum_{k=0}^{n}S_{1}(n,k)x^{k},\quad (\mathrm{see}\ [1-27]),\label{2}	
\end{equation}
where $(x)_{0}=1,\ (x)_{n}=x(x-1)(x-2)\cdots(x-n+1),\ (n\ge 1)$. \\
The Stirling numbers of the second kind are defined by
\begin{equation}
x^{n}=\sum_{k=0}^{n}{n \brace k}(x)_{k},\quad (n\ge 0),\quad (\mathrm{see}\ [5,10]). \label{3}
\end{equation}
Recently, the degenerate Stirling numbers of the second kind are defined by
\begin{equation}
(x)_{n,\lambda}=\sum_{k=0}^{n}{n\brace k}_{\lambda}(x)_{k},\quad (n\ge 0),\quad (\mathrm{see}\ [9,10,19]). \label{3-1}
\end{equation}
Thus, by \eqref{3-1}, we get
\begin{equation}
\frac{1}{k!}\Big(e_{\lambda}(t)-1\Big)^{k}=\sum_{n=k}^{\infty}{n\brace k}_{\lambda}\frac{t^{n}}{n!},\quad (k\ge 0),\quad (\mathrm{see}\ [10]).\label{4}
\end{equation} \par
The degenerate Bell polynomials are defined by
\begin{equation}
\phi_{n,\lambda}(x)=\sum_{k=0}^{n}{n \brace k}_{\lambda}x^{k},\quad (\mathrm{see}\ [9,10,13,17]).\label{5}
\end{equation}
For any $\lambda\in\mathbb{R}$ and $k,n\in\mathbb{N}\cup\{0\}$, Kim-Kim introduced the degenerate Bernstein polynomials which are given by
\begin{equation}
B_{k,n}(x|\lambda)=\binom{n}{k}(x)_{k,\lambda}(1-x)_{n-k,\lambda},\quad (n\ge k\ge 0),\quad (x\in[0,1]).\label{6}
\end{equation}
From \eqref{6}, we note that
\begin{equation}
\frac{1}{k!}(x)_{k,\lambda}t^{k}e_{\lambda}^{1-x}(t)=\sum_{n=k}^{\infty}B_{k,n}(x|\lambda)\frac{t^{n}}{n!},\quad (k\ge 0),\quad (\mathrm{see}\ [16,20]). \label{7}
\end{equation} \par
Assume that $Y$ is a random variable such that the moment generating function of $Y$,
\begin{equation}
E\Big[e^{tY}\Big]=\sum_{n=0}^{\infty}\frac{t^{n}}{n!}E\Big[Y^{n}\Big],\quad (|t|<r),\label{7-1}
\end{equation}
exists for some $r>0$. Let $(Y_{j})_{j\ge 1}$ be a sequence of mutually independent copies of the random variable $Y$, and let
\begin{displaymath}
S_{k}=Y_{1}+Y_{2}+\cdots+Y_{k},\quad (k\ge 1),\quad \mathrm{with}\ S_{0}=0.
\end{displaymath} \par
The probabilistic Bernoulli polynomials associated with $Y$ are  given by
\begin{equation}
\frac{t}{E\big[e^{Yt}\big]-1}\Big(E\big[e^{Yt}\big]\Big)^{x}=\sum_{n=0}^{\infty}B_{n}^{Y}(x)\frac{t^{n}}{n!},\quad (\mathrm{see}\ [13]). \label{8}
\end{equation}
When $x=0$, $B_{n}^{Y}=B_{n}^{Y}(0)$ are called the probabilistic Bernoulli numbers. \\
From \eqref{8}, we have
\begin{equation}
B_{n}^{Y}(x)-B_{n}^{Y}=n\sum_{k=0}^{n-1}\frac{(x)_{k+1}}{k+1}{n-1 \brace k}_{Y}, \label{9}
\end{equation}
where ${n \brace k}_{Y}$ are the probabilistic Stirling numbers of the second kind given by
\begin{equation*}
	\frac{1}{k!}\Big(E\big[e^{Yt}\big]-1\Big)^{k}=\sum_{n=k}^{\infty}{n\brace k}_{Y}\frac{t^{n}}{n!},\quad (\mathrm{see}\ [13,14,17,22]).
\end{equation*}
In particular, $n,k\ge 0$, we have
\begin{equation}
\sum_{k=0}^{n}E\big[S_{k}^{m}\big]=\frac{B_{m+1}^{Y}(n+1)-B_{m+1}^{Y}}{m+1},\quad (\mathrm{see}\ [13]). \label{10}
\end{equation} \par
The probabilistic degenerate Stirling numbers of the second kind associated with $Y$ are introduced by Kim-Kim as
\begin{equation}
\frac{1}{k!}\Big(E\big[e_{\lambda}^{Y}(t)\big]-1\Big)^{k}=\sum_{n=k}^{\infty}{n \brace k}_{Y,\lambda}\frac{t^{n}}{n!},\quad (\mathrm{see}\ [14,17]).\label{11}
\end{equation}
Note that
\begin{displaymath}
{n \brace k}_{1,\lambda}={n \brace k}_{\lambda},\quad (n\ge k\ge 0).
\end{displaymath} \par
In light of \eqref{5}, the probabilistic degenerate Bell polynomials associated with $Y$ are defined by
\begin{equation}
\sum_{n=0}^{\infty}\phi_{n,\lambda}^{Y}(x)\frac{t^{n}}{n!}=e^{x\big(E[e_{\lambda}^{Y}(t)]-1\big)},\quad (\mathrm{see}\ [14,17]). \label{12}	
\end{equation}
The probabilistic degenerate Bernoulli polynomials associated with $Y$ are given by
\begin{equation}
\frac{t}{E\big[e_{\lambda}^{Y}(t)\big]-1}\Big(E\big[e_{\lambda}^{Y}(t)\big]\Big)^{x}=\sum_{n=0}^{\infty}\beta_{n,\lambda}^{Y}(x)\frac{t^{n}}{n!},\quad (\mathrm{see}\ [23]). \label{13}
\end{equation}
When $x=0$, $\beta_{n,\lambda}^{Y}=\beta_{n,\lambda}^{Y}(0)$ are called the probabilistic degenerate Bernoulli numbers associated with $Y$. \\
Also, the probabilistic degenerate Euler polynomials associated with $Y$ are defined by
\begin{equation}
\frac{2}{E\big[e_{\lambda}^{Y}(t)\big]+1}\Big(E\big[e_{\lambda}^{Y}(t)\big]\Big)^{x}=\sum_{n=0}^{\infty}\mathcal{E}_{n,\lambda}^{Y}(x)\frac{t^{n}}{n!},\quad (\mathrm{see}\ [2,3]).\label{14}
\end{equation}
When $x=0$, $\mathcal{E}_{n,\lambda}^{Y}=\mathcal{E}_{n,\lambda}^{Y}(0)$ are called the probabilistic degenerate Euler numbers associated with $Y$. \par
From \eqref{13} and \eqref{14}, we note that
\begin{equation}
\sum_{k=0}^{n}E\big[(S_{k})_{m,\lambda}\big]=\frac{1}{m+1}\Big(\beta_{m+1,\lambda}^{Y}(n+1)-\beta_{m+1,\lambda}^{Y}\Big),\quad (n\ge 0),\label{15}
\end{equation}
and
\begin{equation*}
\sum_{k=0}^{n}(-1)^{k}E\Big[(S_{k})_{m,\lambda}\Big]=\frac{\mathcal{E}_{m,\lambda}^{Y}+\mathcal{E}_{m,\lambda}^{Y}(n+1)}{2},\quad (\mathrm{see}\ [23]).
\end{equation*}
where $n\in\mathbb{N}$ with $n\equiv 0\ (\mathrm{mod}\ 2)$. \\
In addition, $n\in\mathbb{N}$, we have
\begin{equation}
\frac{\beta_{n,\lambda}^{Y}(x)-\beta_{n,\lambda}^{Y}}{n}=\sum_{k=0}^{n-1}\frac{(x)_{k}}{k+1}{n-1 \brace k}_{Y,\lambda},\quad (\mathrm{see}\ [23]).\label{16}
\end{equation} \par
In [7], the probabilistic Bernstein polynomials associated with random variable $Y$ are defined by
\begin{equation}
\sum_{n=k}^{\infty}B_{k,n}^{Y}(x)\frac{t^n}{n!}=\frac{(tx)^{k}}{k!}\Big(E\big[e^{Yt}\big]\Big)^{1-x},\quad (n\ge k\ge 0,\quad x\in[0,1]). \label{17}	
\end{equation}

\section{Probabilistic degenerate Bernstein polynomials associated with random variables}
As in Section 1, we assume that $Y$ is a random variable whose moment generating function exists in a neighborhood of the origin (see \eqref{7-1}). \par
Let us consider the probabilistic degenerate Bernstein polynomials associated with $Y$ which are given by
\begin{equation}
\frac{1}{k!}(x)_{k,\lambda}t^{k}\Big(E\big[e_{\lambda}^{Y}(t)\big]\Big)^{1-x}=\sum_{n=k}^{\infty}B_{k,n}^{Y}(x|\lambda)\frac{t^{n}}{n!},\quad (k\ge 0). \label{18}
\end{equation}
When $Y=1$, $B_{k,\lambda}^{Y}(x|\lambda)=B_{k,n}(x|\lambda),\ (n\ge k\ge 0)$. Note that
\begin{displaymath}
\lim_{\lambda\rightarrow 0}B_{k,n}^{Y}(x|\lambda)=B_{k,n}^{Y}(x),\quad (n\ge k\ge 0).
\end{displaymath} \par
From \eqref{18}, we note that
\begin{align}
\sum_{n=k}^{\infty}B_{k,n}^{Y}(x|\lambda)\frac{t^{n}}{n!}&=\frac{t^{k}}{k!}(x)_{k,\lambda}\Big(E\big[e_{\lambda}^{Y}(t)\big]-1+1\Big)^{1-x} \label{19}\\
&=\frac{(x)_{k,\lambda}}{k!}t^{k}\sum_{j=0}^{\infty}\binom{1-x}{j}\frac{j!}{j!}\Big(E\big[e_{\lambda}^{Y}(t)\big]-1\Big)^{j}\nonumber \\
&=\frac{(x)_{k,\lambda}}{k!}t^{k}\sum_{j=0}^{\infty}\binom{1-x}{j}j!\sum_{n=j}^{\infty}{n \brace j}_{Y,\lambda}\frac{t^{n}}{n!}\nonumber\\
&=\frac{(x)_{k,\lambda}}{k!}t^{k}\sum_{n=0}^{\infty}\bigg(\sum_{j=0}^{n}\binom{1-x}{j}j!{n \brace j}_{Y,\lambda}\bigg)\frac{t^{n}}{n!} \nonumber \\
&=\sum_{n=k}^{\infty}\bigg(\sum_{j=0}^{n-k}\binom{n}{k}\binom{1-x}{j}j!(x)_{k,\lambda}{n-k \brace j}_{Y,\lambda}\bigg)\frac{t^{n}}{n!}. \nonumber
\end{align}
Therefore, by comparing the coefficients on both sides, we obtain the following theorem.
\begin{theorem}
For $n\ge k\ge 0$, we have
\begin{displaymath}
B_{k,n}^{Y}(x|\lambda)=\sum_{j=0}^{n-k}\binom{n}{k}(x)_{k,\lambda}\binom{1-x}{j}j!{n-k \brace j}_{Y,\lambda}.
\end{displaymath}
Note that
\begin{displaymath}
B_{k,n}^{Y}(x)=\lim_{\lambda\rightarrow 0}B_{k,n}^{Y}(x|\lambda)=\sum_{j=0}^{n-k}\binom{n}{k}x^{k}(1-x)_{j}{n-k \brace j}_{Y},\quad (\mathrm{see}\ [7]).
\end{displaymath}
\end{theorem}
For $r\in\mathbb{N}\cup\{0\}$, the probabilistic degenerate Bernoulli polynomials of order $r$ are defined by
\begin{equation}
\bigg(\frac{t}{E\big[e_{\lambda}^{Y}(t)\big]-1}\bigg)^{r}\Big(E\big[e_{\lambda}^{Y}(t)\big]\Big)^{x}=\sum_{n=0}^{\infty}\beta_{n,\lambda}^{(r,Y)}(x)\frac{t^{n}}{n!}.\label{20}
\end{equation}
When $x=0$, $\beta_{n,\lambda}^{(r,Y)}=\beta_{n,\lambda}^{(r,Y)}(0)$ are called the probabilistic degenerate Bernoulli numbers of order $r$. From \eqref{18} and \eqref{20}, we note that
\begin{align}
\sum_{n=k}^{\infty}B_{k,n}^{Y}(x|\lambda)\frac{t^{n}}{n!}&=\frac{t^{k}}{k!}(x)_{k,\lambda}\Big(E\big[e_{\lambda}^{Y}(t)\big]\Big)^{1-x} \label{21} \\
&=\frac{(x)_{k,\lambda}}{k!}\Big(E\big[e_{\lambda}^{Y}(t)\big]-1\Big)^{k}\bigg(\frac{t}{E\big[e_{\lambda}^{Y}(t)\big]-1}\bigg)^{k}\Big(E\big[e_{\lambda}^{Y}(t)\big]\Big)^{1-x}\nonumber\\
&=(x)_{k,\lambda}\sum_{m=k}^{\infty}{m \brace k}_{Y,\lambda}\frac{t^{m}}{m!}\sum_{l=0}^{\infty}\beta_{l,\lambda}^{(k,Y)}(1-x)\frac{t^{l}}{l!}\nonumber\\
&=(x)_{k,\lambda}\sum_{n=k}^{\infty}\bigg(\sum_{m=k}^{n}\binom{n}{m}{m \brace k}_{Y,\lambda}\beta_{n-m,\lambda}^{(k,Y)}(1-x)\bigg)\frac{t^{n}}{n!}. \nonumber
\end{align}
By comparing the coefficients on both sides of \eqref{21}, we obtain the following theorem.
\begin{theorem}
For $n\ge k\ge 0$, we have
\begin{displaymath}
B_{k,n}^{Y}(x|\lambda)=(x)_{k,\lambda}\sum_{m=k}^{n}\binom{n}{m}{m \brace k}_{Y,\lambda}\beta_{n-m,\lambda}^{(k,Y)}(1-x).
\end{displaymath}
\end{theorem}
From \eqref{18}, we note that
\begin{align}
(x)_{k,\lambda}\Big(E\big[e_{\lambda}^{Y}(t)\big]\Big)^{1-x}&=\frac{k!}{t^{k}}\sum_{n=k}^{\infty}B_{k,n}^{Y}(x|\lambda)\frac{t^{n}}{n!} \label{22} \\
&=\frac{k!}{t^{k}}\sum_{n=0}^{\infty}B_{k,n+k}^{Y}(x|\lambda)\frac{t^{n+k}}{(n+k)!}\nonumber\\
&=\sum_{n=0}^{\infty}\frac{B_{k,n+k}^{Y}(x|\lambda)}{\binom{n+k}{n}}\frac{t^{n}}{n!}.\nonumber
\end{align}
On the other hand, by \eqref{11}, we get
\begin{align}
(x)_{k,\lambda}\Big(E\big[e_{\lambda}^{Y}(t)\big]\Big)^{1-x}&=(x)_{k,\lambda}\Big(E\big[e_{\lambda}^{Y}(t)\big]-1+1\Big)^{1-x} \label{23}\\
&=\sum_{j=0}^{\infty}\binom{1-x}{j}j!\frac{1}{j!}\Big(E\big[e_{\lambda}^{Y}(t)\big]-1\Big)^{j} (x)_{k,\lambda}\nonumber\\
&=	\sum_{j=0}^{\infty}\binom{1-x}{j}j!\sum_{n=j}^{\infty}{n\brace j}_{Y,\lambda}\frac{t^{n}}{n!}(x)_{k,\lambda}\nonumber\\
&=\sum_{n=0}^{\infty}\bigg((x)_{k,\lambda} \sum_{j=0}^{n}\binom{1-x}{j}j!{n \brace j}_{Y,\lambda}\bigg)\frac{t^{n}}{n!}.\nonumber
\end{align}
Therefore, by \eqref{22} and \eqref{23}, we obtain the following theorem.
\begin{theorem}
For $n\ge k\ge 0$, we have
\begin{displaymath}
(x)_{k,\lambda}\sum_{j=0}^{n}\binom{1-x}{j}j!{n\brace j}_{Y,\lambda}=\frac{B_{k,n+k}^{Y}(x|\lambda)}{\binom{n+k}{k}}.
\end{displaymath}
\end{theorem}
By \eqref{18}, we get
\begin{align}
&\sum_{n=k}^{\infty}E\Big[(Y)_{n-k,\lambda}\Big]\binom{n}{k}\frac{t^{n}}{n!}=\frac{t^{k}}{k!}E\Big[e_{\lambda}^{Y}(t)\Big]=\frac{\big(E[e_{\lambda}^{Y}(t)]\big)^{x}}{(x)_{k,\lambda}}\sum_{m=k}^{\infty}B_{k,m}^{Y}(x|\lambda)\frac{t^{m}}{m!}	\label{24}\\
&=\sum_{j=0}^{\infty}\frac{1}{(x)_{k,\lambda}}\binom{x}{j}j!\sum_{l=j}^{\infty}{l \brace j}_{Y,\lambda}\frac{t^{l}}{l!}\sum_{m=k}^{\infty}B_{k,m}^{Y}(x|\lambda)\frac{t^{m}}{m!}\nonumber \\
&=\sum_{l=0}^{\infty}\sum_{j=0}^{l}\frac{1}{(x)_{k,\lambda}}\binom{x}{j}j!{l \brace j}_{Y,\lambda}\frac{t^{l}}{l!}\sum_{m=k}^{\infty}B_{k,m}^{Y}(x|\lambda)\frac{t^{m}}{m!}\nonumber \\
&=\sum_{n=k}^{\infty}\bigg(\sum_{l=0}^{n-k}\sum_{j=0}^{l}\frac{1}{(x)_{k,\lambda}}\binom{x}{j}j!{l \brace j}_{Y,\lambda}\binom{n}{l}B_{k,n-l}^{Y}(x|\lambda)\bigg)\frac{t^{n}}{n!}\nonumber\\
&=\sum_{n=k}^{\infty}\bigg(\sum_{l=0}^{n-k}\sum_{j=0}^{l}\frac{\binom{x}{j}}{\binom{x}{k}_{\lambda}}\frac{j!}{k!}{l \brace j}_{Y,\lambda}\binom{n}{l}B_{k,n-l}^{Y}(x|\lambda)\bigg)\frac{t^{n}}{n!}, \nonumber
\end{align}
where $\binom{x}{k}_{\lambda}$ are the degenerate binomial coefficients defined by $\binom{x}{k}_{\lambda}=\frac{(x)_{k,\lambda}}{k!},\ (k\ge 0)$, (see [19]). Therefore, by \eqref{24}, we obtain the following theorem.
\begin{theorem}
For $n\ge k\ge 0$, we have
\begin{displaymath}
E\Big[(Y)_{n-k,\lambda}\Big]= \sum_{l=0}^{n-k}\sum_{j=0}^{l}\frac{\binom{x}{j}}{\binom{x}{k}_{\lambda}}\frac{j!}{k!}{l \brace j}_{Y,\lambda}\frac{\binom{n}{l}}{\binom{n}{k}}B_{k,n-l}^{Y}(x|\lambda).
\end{displaymath}
\end{theorem}
From \eqref{18}, we note that
\begin{align}
&\frac{2}{E\big[e_{\lambda}^{Y}(t)\big]+1}\Big(E\big[e_{\lambda}^{Y}(t)\big]\Big)^{1-x}\frac{t^{k}}{k!}(x)_{k,\lambda}=\frac{2}{E\big[e_{\lambda}^{Y}(t)\big]+1}\sum_{l=k}^{\infty}B_{k,l}^{Y}(x|\lambda)\frac{t^{l}}{l!}\label{29} \\
&=\sum_{m=0}^{\infty}\mathcal{E}_{m,\lambda}^{Y}\frac{t^{m}}{m!}\sum_{l=k}^{\infty}B_{k,l}^{Y}(x|\lambda)\frac{t^{l}}{l!}=\sum_{n=k}^{\infty}\bigg(\sum_{l=k}^{n}\binom{n}{l}\mathcal{E}_{n-l,\lambda}^{Y}B_{k,l}^{Y}(x|\lambda)\bigg)\frac{t^{n}}{n!}.\nonumber
\end{align}
On the other hand, by \eqref{14}, we get
\begin{align}
&\frac{2}{E\big[e_{\lambda}^{Y}(t)\big]+1}\Big(E\big[e_{\lambda}^{Y}(t)\big]\Big)^{1-x}\frac{t^{k}}{k!}(x)_{k,\lambda}=\sum_{n=0}^{\infty}\mathcal{E}_{n,\lambda}^{Y}(1-x)\frac{t^{n}}{n!}\frac{t^{k}}{k!}(x)_{k,\lambda}\label{30}\\
&=\sum_{n=k}^{\infty}(x)_{k,\lambda}\mathcal{E}_{n-k,\lambda}^{Y}(1-x)\binom{n}{k}\frac{t^{n}}{n!}. \nonumber
\end{align}
Therefore, by \eqref{29} and \eqref{30}, we obtain the following theorem.
\begin{theorem}
For $n\ge k\ge 0$, we have
\begin{displaymath}
\mathcal{E}_{n-k,\lambda}^{Y}(1-x)=\sum_{l=k}^{n}\frac{\binom{n}{l}}{\binom{n}{k}\binom{x}{k}_{\lambda}}\mathcal{E}_{n-l,\lambda}^{Y}\frac{B_{k,l}^{Y}(x|\lambda)}{k!}.
\end{displaymath}
\end{theorem}
From \eqref{18} and \eqref{20}, we have
\begin{align}
&\sum_{n=k}^{\infty}B_{k,n}^{Y}(x|\lambda)\frac{t^{n}}{n!}=\frac{1}{k!}t^{k}(x)_{k,\lambda}\Big(E\big[e_{\lambda}^{Y}(t)\big]\Big)^{1-x} \label{31} \\
&=(x)_{k,\lambda}\bigg(\frac{t}{E\big[e_{\lambda}^{Y}(t)\big]-1}\bigg)^{k}\Big(E\big[e_{\lambda}^{Y}(t)\big]\Big)^{kx}\frac{1}{k!}\Big(E\big[e_{\lambda}^{Y}(t)\big]-1\Big)^{k}\Big(E\big[e_{\lambda}^{Y}(t)\big]\Big)^{1-x-2k}\nonumber\\
&=(x)_{k,\lambda}\sum_{l=0}^{\infty}\beta_{l,\lambda}^{(k,Y)}(kx)\frac{t^{l}}{l!}\frac{1}{k!}\sum_{j=0}^{\infty}\binom{1-x-kx}{j}j!\frac{1}{j!}\Big(E\big[e_{\lambda}^{Y}(t)-1\Big)^{k+j}\nonumber 	\\
&=(x)_{k,\lambda}\sum_{l=0}^{\infty}\beta_{l,\lambda}^{(k,Y)}(kx)\frac{t^{l}}{l!}\sum_{j=0}^{\infty}\binom{k+j}{j}\binom{1-x-kx}{j}j!\sum_{m=k+j}^{\infty}{m \brace k+j}_{Y,\lambda}\frac{t^{m}}{m!}\nonumber \\
&=(x)_{k,\lambda}\sum_{l=0}^{\infty}\beta_{l,\lambda}^{(k,Y)}(kx)\frac{t^{l}}{l!}\sum_{m=k}^{\infty}\sum_{j=0}^{m-k}\binom{k+j}{j}\binom{1-x-kx}{j}j!{m \brace k+j}_{Y,\lambda}\frac{t^{m}}{m!}\nonumber \\
&=\sum_{n=k}^{\infty}\bigg((x)_{k,\lambda}\sum_{m=k}^{n}\sum_{j=0}^{m-k}\binom{n}{m}\beta_{n-m,\lambda}^{(k,Y)}(kx)\binom{k+j}{j}\binom{1-x-kx}{j}j!{m \brace k+1}_{Y,\lambda}\bigg)\frac{t^{n}}{n!}.\nonumber
\end{align}
Therefore, by comparing the coefficients on both sides of \eqref{31}, we obtain the following theorem.
\begin{theorem}
For $n\ge k\ge 0$, we have
\begin{displaymath}
B_{k,n}^{Y}(x|\lambda)= (x)_{k,\lambda}\sum_{m=k}^{n}\sum_{j=0}^{m-k}\binom{n}{m}\beta_{n-m,\lambda}^{(k,Y)}(kx)\binom{k+j}{j}\binom{1-x-kx}{j}{m \brace k+1}_{Y,\lambda}j!.
\end{displaymath}
\end{theorem}
Now, we observe that
\begin{align}
&\sum_{n=k}^{\infty}B_{k,n}^{Y}(x|\lambda)\frac{t^{n}}{n!}=\frac{t^{k}}{k!}(x)_{k,\lambda}\Big(E\big[e_{\lambda}^{Y}(t)\big]\Big)^{1-x}=\frac{t^{k}}{k!}(x)_{k,\lambda}e^{(1-x)\log E[e_{\lambda}^{Y}(t)]}\label{32}\\
&=\frac{t^{k}}{k!}(x)_{k,\lambda}\sum_{l=0}^{\infty}(1-x)^{l}\frac{1}{l!}\Big(\log\big(E[e_{\lambda}^{Y}(t)]-1+1\big)\Big)^{l}\nonumber\\
&=\frac{t^{k}}{k!}(x)_{k,\lambda}\sum_{l=0}^{\infty}(1-x)^{l}\sum_{j=l}^{\infty}S_{1}(j,l)\frac{1}{j!}\Big(E\big[e_{\lambda}^{Y}(t)]-1\Big)^{j} \nonumber \\
&=\frac{t^{k}}{k!}(x)_{k,\lambda}\sum_{j=0}^{\infty}\bigg(\sum_{l=0}^{j}(1-x)^{l}S_{1}(j,l)\bigg)\sum_{n=j}^{\infty}{n \brace j}_{Y,\lambda}\frac{t^{n}}{n!}\nonumber\\
&=\sum_{n=k}^{\infty}\bigg((x)_{k,\lambda}\sum_{j=0}^{n-k}\sum_{l=0}^{j}(1-x)^{l}S_{1}(j,l){n-k \brace j}_{Y,\lambda}\binom{n}{k}\bigg)\frac{t^{n}}{n!}.\nonumber
\end{align}
Therefore, by \eqref{32}, we obtain the following theorem.
\begin{theorem}
For $n\ge k\ge 0$, we have
\begin{displaymath}
\frac{1}{k!}B_{k,n}^{Y}(x|\lambda)=\binom{x}{k}_{\lambda}\binom{n}{k}\sum_{j=0}^{n-k}\sum_{l=0}^{j}(1-x)^{l}S_{1}(j,l){n-k \brace j}_{Y,\lambda}.
\end{displaymath}
\end{theorem}
Let $Y$ be the Poisson random variable with parameter $\alpha$. Then we have
\begin{equation}
E\Big[e_{\lambda}^{Y}(t)\Big]=\sum_{n=0}^{\infty}e_{\lambda}^{n}(t)\frac{\alpha^{n}}{n!}e^{-\alpha}=e^{\alpha(e_{\lambda}(t)-1)}=\sum_{n=0}^{\infty}\phi_{n,\lambda}(\alpha)\frac{t^{n}}{n!}. \label{33}	
\end{equation}
By \eqref{33}, we get
\begin{equation}
\Big(E\big[e_{\lambda}^{Y}(t)\big]\Big)^{1-x}=e^{\alpha(1-x)(e_{\lambda}(t)-1)}=\sum_{n=0}^{\infty}\phi_{n,\lambda}\big(\alpha(1-x)\big)\frac{t^{n}}{n!}.\label{34}
\end{equation}
From \eqref{18} and \eqref{34}, we have
\begin{align}
\sum_{n=k}^{\infty}B_{k,n}^{Y}(x|\lambda)\frac{t^{n}}{n!}&=\frac{t^{k}}{k!}(x)_{k,\lambda}\Big(E\big[e_{\lambda}^{Y}(t)\big]\Big)^{1-x} \label{35}\\
&= \frac{t^{k}}{k!}(x)_{k,\lambda}\sum_{n=0}^{\infty}\phi_{n,\lambda}\big(\alpha(1-x)\big)\frac{t^{n}}{n!}\nonumber\\
&= \frac{t^{k}}{k!}(x)_{k,\lambda} \sum_{n=k}^{\infty}\phi_{n-k,\lambda}\big(\alpha(1-x)\big)\frac{t^{n-k}}{(n-k)!}\nonumber \\
&=\sum_{n=k}^{\infty}(x)_{k,\lambda}\phi_{n-k,\lambda}\big(\alpha(1-x)\big)\binom{n}{k}\frac{t^{n}}{n!}.\nonumber
\end{align}
Therefore, by \eqref{35}, we obtain the following theorem.
\begin{theorem}
Let $Y$ be the Poisson random variable with parameter $\alpha(>0)$. For $n\ge k\ge 0$, we have
\begin{displaymath}
\frac{B_{k,n}^{Y}(x|\lambda)}{k!}=\binom{x}{k}_{\lambda}\binom{n}{k}\phi_{n-k,\lambda}\big(\alpha(1-x)\big).
\end{displaymath}
\end{theorem}
Assume that $Y$ is the Poisson random variable with parameter $\alpha(>0)$. From \eqref{34}, we have
\begin{align}
\sum_{n=k}^{\infty}B_{k,n}^{Y}(x|\lambda)\frac{t^{n}}{n!}&=\frac{t^{k}}{k!}(x)_{k,\lambda}\Big(E\big[e_{\lambda}^{Y}(t)\big]\Big)^{1-x} \label{36} \\
&=\frac{t^{k}}{k!}(x)_{k,\lambda}e^{\alpha(1-x)(e_{\lambda}(t)-1)} \nonumber\\
&= \frac{t^{k}}{k!}(x)_{k,\lambda}\sum_{j=0}^{\infty}\alpha^{j}(1-x)^{j}\frac{1}{j!}\big(e_{\lambda}(t)-1\big)^{j}\nonumber\\
&= \frac{t^{k}}{k!}(x)_{k,\lambda}\sum_{j=0}^{\infty}\alpha^{j}(1-x)^{j}\sum_{n=j}^{\infty}{n \brace j}_{\lambda}\frac{t^{n}}{n!} \nonumber \\
&= \frac{t^{k}}{k!}(x)_{k,\lambda}\sum_{n=0}^{\infty}\bigg(\sum_{j=0}^{n}\alpha^{j}(1-x)^{j}{n \brace j}_{\lambda}\bigg)\frac{t^{n}}{n!} \nonumber\\
&=\sum_{n=k}^{\infty}\bigg((x)_{k,\lambda}\sum_{j=0}^{n-k}\binom{n}{k}\alpha^{j}(1-x)^{j}{n-j \brace j}_{\lambda}\bigg)\frac{t^{n}}{n!}. \nonumber
\end{align}
Therefore, by \eqref{36}, we obtain the following theorem.
\begin{theorem}
Assume that $Y$ is the Poisson random variable with parameter $\alpha(>0)$. For $n\ge k\ge 0$, we have
\begin{displaymath}
\frac{B_{k,n}^{Y}(x|\lambda)}{k!}=\binom{x}{k}_{\lambda}\binom{n}{k}\sum_{j=0}^{n-k}\alpha^{j}(1-x)^{j}{n-k \brace j}_{\lambda}.
\end{displaymath}
\end{theorem}
Let $Y$ be the Bernoulli random variable with probability of success $p$. Then we have
\begin{equation}
E\Big[e_{\lambda}^{Y}(t)\Big]=e_{\lambda}(t)p+1-p=p(e_{\lambda}(t)-1)+1. \label{37}	
\end{equation}
From \eqref{18} and \eqref{37}, we have
\begin{align}
\sum_{n=k}^{\infty}B_{k,n}^{Y}(x|\lambda)\frac{t^{n}}{n!}&=\frac{t^{k}}{k!}(x)_{k,\lambda}\Big(E\big[e_{\lambda}^{Y}(t)\big]\Big)^{1-x} \label{38}\\
&=\frac{t^{k}}{k!}(x)_{k,\lambda}\Big(p(e_{\lambda}(t)-1)+1\Big)^{1-x} \nonumber\\
&=\frac{t^{k}}{k!}(x)_{k,\lambda}\sum_{l=0}^{\infty}\binom{1-x}{l}l!\frac{1}{l!}\big(p(e_{\lambda}(t)-1)\big)^{l} \nonumber\\
&=\frac{t^{k}}{k!}(x)_{k,\lambda}\sum_{l=0}^{\infty}p^{l}\binom{1-x}{l}l!\sum_{n=l}^{\infty}{n \brace l}_{\lambda}\frac{t^{n}}{n!} \nonumber \\
&=\frac{t^{k}}{k!}(x)_{k,\lambda}\sum_{n=0}^{\infty}\bigg(\sum_{l=0}^{n}\binom{1-x}{l}l!p^{l}{n \brace l}_{\lambda}\bigg)\frac{t^{n}}{n!}\nonumber \\
&=\sum_{n=k}^{\infty}\binom{n}{k}(x)_{k,\lambda}\sum_{l=0}^{n-k}\binom{1-x}{l}l!p^{l}{n-k \brace l}_{\lambda}\frac{t^{n}}{n!}.
\end{align}
Therefore, by \eqref{38}, we obtain the following theorem.
\begin{theorem}
Let $Y$ be the Bernoulli random variable with probability of success $p$. For $n\ge k\ge 0$, we have
\begin{displaymath}
\frac{B_{k,n}^{Y}(x|\lambda)}{k!}=\binom{n}{k}\binom{x}{k}_{\lambda}\sum_{l=0}^{n-k}p^{l}\binom{1-x}{l}l!{n-k \brace l}_{\lambda}.
\end{displaymath}
\end{theorem}
Assume that $Y$ is the binomial random variable with parameters $m,p$. Then we have
\begin{align}
E\Big[e_{\lambda}^{Y}(t)\Big]&=\sum_{k=0}^{m}\binom{m}{k}p^{k}(1-p)^{m-k}e_{\lambda}^{k}(t) \label{39}\\
&=\big(pe_{\lambda}(t)+1-p\big)^{m}=\big(p(e_{\lambda}(t)-1)+1\big)^{m}. \nonumber	
\end{align}
By \eqref{18} and \eqref{39}, we get
\begin{align}
\sum_{n=k}^{\infty}B_{k,n}^{Y}(x|\lambda)\frac{t^{n}}{n!}&=\frac{t^{k}}{k!}(x)_{k,\lambda}\Big(E\big[e_{\lambda}^{Y}(t)\big]\Big)^{1-x}\label{40}	\\
&=\frac{t^{k}}{k!}(x)_{k,\lambda}\Big(p(e_{\lambda}(t)-1)+1\Big)^{(1-x)m} \nonumber\\
&=\frac{t^{k}}{k!}(x)_{k,\lambda}\sum_{j=0}^{\infty}\binom{m(1-x)}{j}j!p^{j}\frac{1}{j!}\big(e_{\lambda}(t)-1\big)^{j} \nonumber\\
&=\frac{t^{k}}{k!}(x)_{k,\lambda}\sum_{j=0}^{\infty}\binom{m(1-x)}{j}j!p^{j}\sum_{n=j}^{\infty}{n\brace j}_{\lambda}\frac{t^{n}}{n!}\nonumber\\
&=\frac{t^{k}}{k!}(x)_{k,\lambda}\sum_{n=0}^{\infty}\bigg(\sum_{j=0}^{n}\binom{m(1-x)}{j}j!p^{j}{n \brace j}_{\lambda}\bigg)\frac{t^{n}}{n!}\nonumber\\
&=\sum_{n=k}^{\infty}\binom{n}{k}(x)_{k,\lambda}\sum_{j=0}^{n-k}\binom{m(1-x)}{j}j!p^{j}{n-j \brace j}_{\lambda}\frac{t^{n}}{n!}.\nonumber
\end{align}
Therefore, by \eqref{40}, we obtain the following theorem.
\begin{theorem}
Assume that $Y$ is the binomial random variable with parameters $m,p(>0)$. For $n\ge k\ge 0$, we have
	\begin{displaymath}
		\frac{B_{k,n}^{Y}(x|\lambda)}{k!}=\binom{n}{k}\binom{x}{k}_{\lambda}\sum_{j=0}^{n-k}\binom{m(1-x)}{j}j!p^{j}{n-k \brace j}_{\lambda}.
	\end{displaymath}
\end{theorem}
\section{Conclusion}
Both degenerate Bernstein polynomials (see [15,16,20]) and probabilistic Bernstein polynomials (see [7]) were studied earlier. Let $Y$ be a random variable whose moment generating function exists in a neighborhood of the origin (see \eqref{7-1}). In this paper, we introduced the probabilistic degenerate Bernstein polynomials $B_{k,n}^{Y}(x|\lambda)$ which can be viewed as both probabilistic extension of the degenerate Bernstein polynomials and degenerate version of the probabilistic Bernstein polynomials. We derived several explicit expressions and related identities for those polynomials in connection with other polynomials, including the probabilistic degenerate Bernoulli polynomials associated with $Y$ and the probabilistic degenerate Euler polynomials associated with $Y$. \par
As our future research projects, we would like to continue to explore degenerate versions and probabilistic extensions of many special numbers and polynomials.

\end{document}